\documentclass[reqno,onecolumn]{amsproc}%
\usepackage{amsfonts}
\usepackage{amsmath}
\usepackage{amssymb}
\usepackage{graphicx}
\usepackage[none]{hyphenat}%
\theoremstyle{plain}

\numberwithin{equation}{section}
\begin{document}
\title{}

\begin{center}
{\Huge Details for "Least-Squares Prices of Games"}

\bigskip

{\large Yukio Hirashita}

\bigskip

\bigskip

\textbf{Abstract}

\bigskip
\end{center}

\noindent This paper is intended to help the readers to understand the article:

\ \ Y. Hirashita, Least-Squares Prices of Games, Preprint,

\ \ arXiv:math.OC/0703079 (2007).

\bigskip



\bigskip

\bigskip

\noindent\textbf{Remark to Theorem 1.1. }To understand Theorem 1.1, it is
useful to run the following Mathematica program. For any positive values of a,
b, and r, the theoretical growth rate E\symbol{94}r and the simulated growth
rate (\textquotedblleft geometric mean\textquotedblright) x\symbol{94}%
(1/Repeat) are almost equal.
\begin{verbatim}
       a=19;b=1;
       r=0.05; Print["theoretical growth rate = ",E^r];
       EA=(a+b)/2;k=(1-Sqrt[1-1/E^(2r)])/2;
       If[a<b,c=a;a=b;b=c];uA=Sqrt[a*b]/E^r;t=1;
       If[EA>Sqrt[a*b]*E^r,uA=k*a+(1-k)*b;t=uA(EA-uA)/((a-uA)(uA-b))];
       x=1;Repeat=100000;
       Do[If[Random[]<0.5,x=x*t*a/uA+x*(1-t),x=x*t*b/uA+x*(1-t)],{n,1,Repeat}];
       Print["simulated growth rate = ",x^(1/Repeat)];
\end{verbatim}

\bigskip

\noindent\textbf{Proof of Theorem 1.1.} From Remark 3.1, in the case where
\begin{align*}
\frac{1}{e^{r}}\exp(\int\log a(x)dF(x))  &  =\frac{1}{e^{r}}\exp(\frac{\log
a}{2}+\frac{\log b}{2})=\frac{\sqrt{ab}}{e^{r}}\\
&  \leq\frac{1}{\int\frac{1}{a(x)}dF(x)}=\frac{1}{\frac{1}{2a}+\frac{1}{2b}%
}=\frac{2ab}{a+b},
\end{align*}
that is, in the case where $(a+b)/(2\sqrt{ab})=E/\sqrt{ab}\leq e^{r}$, the
price is given by $\exp(\int\log a(x)dF(x))/e^{r}=\sqrt{ab}/e^{r}$, and the
optimal proportion of investment is $1$. Otherwise, the price $u>0$ and the
optimal proportion of investment $t_{u}>0$ are determined by the simultaneous
equations
\[
\left\{
\begin{array}
[c]{c}%
\exp(\int\log(\frac{a(x)t_{u}}{u}-t_{u}+1)dF(x))=\frac{\sqrt{(at_{u}%
-ut_{u}+u)(bt_{u}-ut_{u}+u)}}{u}=e^{r},\\
\int\frac{a(x)-u}{a(x)t_{u}-ut_{u}+u}dF(x)=\frac{a-u}{2(at_{u}-ut_{u}%
+u)}+\frac{b-u}{2(bt_{u}-ut_{u}+u)}=0.
\end{array}
\right.
\]
It is not difficult to verify that the solutions are given by $u=\kappa
a+(1-\kappa)b$ and $t_{u}=u(E-u)/((a-u)(u-b))$, where $\kappa$ $:=(1-\sqrt
{1-1/e^{2r}})/2$.\hfill$\square$

\bigskip

\noindent\textbf{Lemma D.1. }$u_{r}^{kA}=ku_{r}^{A}$\textit{ for }$k>0$.

\noindent\textit{Proof.} From Remark 3.1, in the case where $\exp(\int\log
a(x)dF(x))/e^{r}\leq1/\int1/a(x)dF(x)$, we have $\exp(\int\log
(ka(x))dF(x))/e^{r}$ $\leq1/\int1/(ka(x))dF(x)$ and $u_{r}^{kA}=\exp(\int
\log(ka(x))$ $dF(x))/e^{r}$ $=k\exp(\int\log a(x)dF)/e^{r}$ $=ku_{r}^{A}.$ In
the other case, as $a(x)/u=ka(x)/(ku)$ and $(a(x)-u)$ $/(a(x)t_{u}-ut_{u}+u)$
$=(ka(x)-ku)/(ka(x)t_{u}-kut_{u}+ku)$, the pattern of the simultaneous
equations remain unchanged.\hfill$\square$

\bigskip

\noindent\textbf{Lemma D.2.} $T$\textit{ is convex}.

\noindent\textit{Proof}. By definition, for each $(t_{i})$ and $(t^{\prime}%
{}_{i})$ in $T$ we have%
\[
\frac{u_{r}^{\sum_{i=1}^{n}p_{i}A_{i}}}{\sum_{i=1}^{n}p_{i}(u_{r}^{A_{i}%
}+t_{i}(E_{A_{i}}/e^{r}-u_{r}^{A_{i}}))}\leq1\text{ and }\frac{u_{r}%
^{\sum_{i=1}^{n}p_{i}A_{i}}}{\sum_{i=1}^{n}p_{i}(u_{r}^{A_{i}}+t^{\prime}%
{}_{i}(E_{A_{i}}/e^{r}-u_{r}^{A_{i}}))}\leq1
\]
for each $(p_{i})\in Q$. Thus,
\[
\frac{u_{r}^{\sum_{i=1}^{n}p_{i}A_{i}}}{\sum_{i=1}^{n}p_{i}(u_{r}^{A_{i}%
}+(qt_{i}+(1-q)t^{\prime}{}_{i})(E_{A_{i}}/e^{r}-u_{r}^{A_{i}}))}\leq1\text{
\ \ }(0\leq q\leq1),
\]
which implies the conclusion.\hfill$\square$

\bigskip

\noindent\textbf{Lemma D.3.} $T$\textit{ is closed}.

\noindent\textit{Proof.} Put%
\[
f_{(p_{i})}((t_{i})):=\frac{u_{r}^{\sum_{i=1}^{n}p_{i}A_{i}}}{\sum_{i=1}%
^{n}p_{i}(u_{r}^{A_{i}}+t_{i}(E^{A_{i}}/e^{r}-u_{r}^{A_{i}}))},
\]
then, as $\sum_{i=1}^{n}p_{i}u_{r}^{A_{i}}\geq\min_{1\leq i\leq n}u_{r}%
^{A_{i}}>0$, $f_{(p_{i})}((t_{i}))$ is continuous with respect to $(t_{i})\in
S$. Therefore, $\{(t_{i})\in S$ : $f_{(p_{i})}((t_{i}))\leq1\}$ is closed in
$S$ for each $(p_{i})\in Q$. Thus, $\cap_{(p_{i})\in Q}$ $\{(t_{i})\in S$ :
$f_{(p_{i})}((t_{i}))\leq1\}=\{(t_{i})\in S$ : $L((t_{i}))$ $\leq1\}=T$ is closed.

\ \hfill$\square$

\bigskip

\noindent\textbf{Remark to Definition 2.1. }In Definition 2.1, we can write
\[
L((t_{i})):=\max_{(p_{i})\in Q}\frac{u_{r}^{\sum_{i=1}^{n}p_{i}A_{i}}}%
{\sum_{i=1}^{n}p_{i}(u_{r}^{A_{i}}+t_{i}(E^{A_{i}}/e^{r}-u_{r}^{A_{i}}%
))}\text{ \ \ }((t_{i})\in S),
\]
because $u_{r}^{\sum_{i=1}^{n}p_{i}A_{i}}$ is continuous with respect to
$(p_{i})\in Q$ (see Theorem D.19). Moreover, by Berge's maximum theorem
[\textbf{8}, Theorem 2.1], $L((t_{i}))$ is continuous with respect to
$(t_{i})\in S$.

\bigskip

\noindent\textbf{Lemma D.4. }$L((x_{i}))=1$ \textit{for} $u_{r}^{A_{i},\text{
}\Omega}=u_{r}^{A_{i}}+x_{i}(E^{A_{i}}/e^{r}-u_{r}^{A_{i}})$ $(0\leq i\leq1)$.

\noindent\textit{Proof.} From the continuity of $u_{r}^{\sum_{i=1}^{n}%
p_{i}A_{i}}$, $f_{(p_{i})}((t_{i}))=u_{r}^{\sum_{i=1}^{n}p_{i}A_{i}}%
/\sum_{i=1}^{n}p_{i}(u_{r}^{A_{i}}$ $+t_{i}(E^{A_{i}}/e^{r}-u_{r}^{A_{i}}))$
is uniformly continuous with respect to $((t_{i}),(p_{i}))$ on the compact set
$S\times Q$. Assume $L((x_{i}))<1$ and choose $(q_{i})\in Q$ such that
\[
L((x_{i}))=\max_{(p_{i})\in Q}\frac{u_{r}^{\sum_{i=1}^{n}p_{i}A_{i}}}%
{\sum_{i=1}^{n}p_{i}(u_{r}^{A_{i}}+x_{i}(E^{A_{i}}/e^{r}-u_{r}^{A_{i}}%
))}=\frac{u_{r}^{\sum_{i=1}^{n}q_{i}A_{i}}}{\sum_{i=1}^{n}q_{i}u_{r}%
^{A_{i},\text{ }\Omega}}<1.
\]
If $x_{j}>0$ exists, then there is a $0<\varepsilon<1$ such that
$L((x_{i}^{\prime}))<1$, where $x_{j}^{\prime}=\varepsilon x_{j}$,
$x_{i}^{\prime}=x_{i}$ $(i\neq j)$ and $\sum_{i=1}^{n}\left(  x_{i}^{\prime
}\right)  ^{2}<\sum_{i=1}^{n}x_{i}^{2}$, which is a contradiction. On the
other hand, if $x_{i}=0$ for each $0\leq i\leq1$, then $L((0))<1$, which is
also a contradiction.\hfill$\square$

\bigskip

\noindent\textbf{Notation D.5.} For two games $A=(a(x)$, $dF(x))$ and
$B=(b(x)$, $dF(x))$, we use the following notation:

$f(p):=\exp(\int\log(pa(x)+(1-p)b(x))dF(x))/e^{r}$ $\ \ (0\leq p\leq1).$

$g(p)$ is defined by $u$ of the simultaneous equations:%
\[
\left\{
\begin{array}
[c]{c}%
\exp(\int\log(\frac{pa(x)+(1-p)b(x)}{u}t_{u}-t_{u}+1)dF(x))=e^{r},\\
\int\frac{pa(x)+(1-p)b(x)-u}{\left(  pa(x)+(1-p)b(x)\right)  t_{u}-ut_{u}%
+u}dF(x)=0\text{ \ \ }(0\leq p\leq1).
\end{array}
\right.
\]
$\ \ \ \ \ \ h(p):=1/(\int1/(pa(x)+(1-p)b(x))dF(x))$ $\ \ (0\leq p\leq1).$%
\[
u(p):=u_{r}^{pA+(1-p)B}=\left\{
\begin{array}
[c]{c}%
f(p)\text{ if }f(p)\leq h(p),\\
g(p)\text{ if }f(p)>h(p)\ \ \ (0\leq p\leq1).
\end{array}
\right.
\]

\bigskip

\noindent\textbf{Lemma D.6.} $f(p)\leq g(p)$\textit{, if }$g(p)$
\textit{exists}.

\noindent\textit{Proof}. As the function $\exp(\int\log(\left(
pa(x)+(1-p)b(x)\right)  t/u-t+1)dF(x))$ is concave with respect to $t$ (see
[\textbf{3}, Lemma 4.7]), it reaches its maximum at $t=t_{u}$. Therefore, we
have%
\begin{align*}
e^{r}  &  =\exp(\int\log(\frac{pa(x)+(1-p)b(x)}{g(p)}t_{g(p)}-t_{g(p)}%
+1)dF(x))\\
&  \geq\frac{\exp(\int\log(pa(x)+(1-p)b(x))dF(x))}{g(p)}.
\end{align*}
On the other hand, we have $e^{r}=\exp(\int\log(pa(x)+(1-p)b(x))dF(x))/f(p).$
Therefore, $1/f(p)\geq1/g(p)$, which implies the conclusion.\hfill$\square$

\bigskip

\noindent\textbf{Lemma D.7.} \textit{The following four properties are
equivalent at a point} $p\in\lbrack0,$ $1]$.

(1) $f(p)=g(p)=h(p)$.

(2) $f(p)=g(p)$.

(3) $f(p)=h(p)$.

(4) $g(p)=h(p)$.

\noindent\textit{Proof}. (3) $\Longrightarrow$ (2). We write
$c(x):=pa(x)+(1-p)b(x)$. From $f(p)=h(p)$, we have%

\[
\frac{\exp(\int\log c(x)dF(x))}{e^{r}}=\frac{1}{\int\frac{1}{c(x)}dF(x)}.
\]
Write $u$ for this value and put $t_{u}:=1$. Then, we obtain%
\[
\left\{
\begin{array}
[c]{c}%
\exp(\int\log(\frac{c(x)}{u}t_{u}-t_{u}+1)dF(x))=\exp(\int\log\frac{c(x)}%
{u}dF(x))=e^{r},\\
\int\frac{c(x)-u}{c(x)}dF(x)=1-u\int\frac{1}{c(x)}dF(x)=0.
\end{array}
\right.
\]
Therefore, by the uniqueness of the solutions (see [\textbf{3}, Section 6]),
we have $u=g(p)$.

(4) $\Longrightarrow$ (2). Put $u:=g(p)$ and $H:=1/h(p)$. Then, $u=h(p)$
implies $u=1/H$. From [\textbf{3}, Lemmas 4.12, 4.16, and 4.21], we obtain
$e^{r}=H\exp(\int\log c(x)dF(x))$, which implies $h(p)=f(p)$.

The other cases can be obtained in a similar fashion.\hfill$\square$

\bigskip

\noindent\textbf{Lemma D.8.} $f(p)$\textit{ is concave on} $[0$, $1]$.

\noindent\textit{Proof.}\textbf{ }Let $\{p,$ $q,$ $\lambda\}\subset
\lbrack0,1]$. By the fact that $\lambda\exp(\int\log a(x)dF(x))=\exp(\int
\log(\lambda$ $a(x))dF(x))$ and using [\textbf{2}, Theorem 185], we obtain
\begin{align*}
&  \lambda f(p)+(1-\lambda)f(q)\\
&  =\frac{\left(
\begin{array}
[c]{c}%
\exp(\int\log(\lambda pa(x)+\lambda(1-p)b(x))dF(x))\\
+\exp(\int\log((1-\lambda)qa(x)+(1-\lambda)(1-q)b(x))dF(x))
\end{array}
\right)  }{e^{r}}\\
&  \leq\frac{\exp(\int\log((\lambda p+(1-\lambda)q)a(x)+(1-(\lambda
p+(1-\lambda)q))b(x))dF(x))}{e^{r}}\\
&  =f(\lambda p+(1-\lambda)q).
\end{align*}
Thus, we have the conclusion.\hfill$\square$

\bigskip

\noindent\textbf{Lemma D.9. }$g(p)$\textit{ is concave on }$[0$\textit{, }%
$1]$\textit{ if }$g(p)$\textit{ exists for each} $p\in\lbrack0,$ $1]$.

\noindent\textit{Proof.}\textbf{ }Let $\{p,$ $q,$ $\lambda\}\subset
\lbrack0,1]$, $C:=(c(x)$, $dF(x)),$ and $D:=(d(x)$, $dF(x))$, where $c(x)$
$:=pa(x)+(1-p)b(x)$ and $d(x):=qa(x)+(1-q)b(x)$. Notice that%
\begin{align*}
e^{r}  &  =\exp(\int\log(\frac{c(x)}{u^{C}}t_{C}-t_{C}+1)dF(x))\text{, }\\
e^{r}  &  =\exp(\int\log(\frac{d(x)}{u^{D}}t_{D}-t_{D}+1)dF(x)),
\end{align*}
where $u^{C}:=g(p)$, $u^{D}:=g(q)$, $t_{C}:=t_{u^{C}}$, and $t_{D}:=t_{u^{D}}%
$. Be careful that $u_{r}^{C}\in\{f(p),$ $g(p)\}$ is not necessarily equal to
$u^{C}$. Put $\mu:=\lambda/(\lambda+(1-\lambda)t_{C}u^{D}/(t_{D}u^{C}))$ and
$\widehat{t}$ $:=\mu t_{C}$ $+(1-\mu)t_{D}$. Then, using [\textbf{2}, Theorem
185], we have%
\begin{align*}
e^{r}  &  =\mu\exp(\int\log(\frac{c(x)}{u^{C}}t_{C}-t_{C}+1)dF(x))\\
&  +(1-\mu)\exp(\int\log(\frac{d(x)}{u^{D}}t_{D}-t_{D}+1)dF(x))\\
&  =\exp(\int\log(\frac{\mu c(x)}{u^{C}}t_{C}-\mu t_{C}+\mu)dF(x))\\
&  +\exp(\int\log(\frac{(1-\mu)d(x)}{u^{D}}t_{D}-(1-\mu)t_{D}+(1-\mu))dF(x))\\
&  \leq\exp(\int\log(\frac{\mu c(x)}{u^{C}}t_{C}+\frac{(1-\mu)d(x)}{u^{D}%
}t_{D}-(\mu t_{C}+(1-\mu)t_{D})+1)dF(x))\\
&  =\exp(\int\log(\frac{\lambda c(x)+(1-\lambda)d(x)}{\lambda u^{C}%
+(1-\lambda)u^{D}}\widehat{t}-\widehat{t}+1)dF(x)).
\end{align*}
On the other hand, we have
\begin{align*}
e^{r}  &  =\exp(\int\log(\frac{\lambda c(x)+(1-\lambda)d(x)}{u^{\lambda
C+(1-\lambda)D}}t_{u^{\lambda C+(1-\lambda)D}}-t_{u^{\lambda C+(1-\lambda)D}%
}+1)dF(x))\\
&  \geq\exp(\int\log(\frac{\lambda c(x)+(1-\lambda)d(x)}{u^{\lambda
C+(1-\lambda)D}}\widehat{t}-\widehat{t}+1)dF(x)).
\end{align*}
Therefore,%
\begin{align*}
&  \exp(\int\log(\frac{\lambda c(x)+(1-\lambda)d(x)}{u^{\lambda C+(1-\lambda
)D}}\widehat{t}-\widehat{t}+1)dF(x))\\
&  \leq\exp(\int\log(\frac{\lambda c(x)+(1-\lambda)d(x)}{\lambda
u^{C}+(1-\lambda)u^{D}}\widehat{t}-\widehat{t}+1)dF(x)),
\end{align*}
which implies $u^{\lambda C+(1-\lambda)D}\geq\lambda u^{C}+(1-\lambda)u^{D}$,
and also $g(\lambda p+(1-\lambda)q)\geq\lambda g(p)+(1$ $-\lambda)g(q).$%
\hfill$\square$

\bigskip

\noindent\textbf{Lemma D.10.} $h(p)$\textit{ is concave on} $[0$, $1]$.

\noindent\textit{Proof.}\textbf{ }Let $\{p,$ $q,$ $\lambda\}\subset
\lbrack0,1]$. Then, using [\textbf{2}, Theorem 214], we obtain
\begin{align*}
&  \lambda h(p)+(1-\lambda)h(q)=\frac{\lambda}{\int\frac{1}{pa(x)+(1-p)b(x)}%
dF(x)}+\frac{1-\lambda}{\int\frac{1}{qa(x)+(1-q)b(x)}dF(x)}\\
&  =\frac{1}{\int\frac{1}{\lambda(pa(x)+(1-p)b(x))}dF(x)}+\frac{1}{\int
\frac{1}{(1-\lambda)(qa(x)+(1-q)b(x))}dF(x)}\\
&  \leq\frac{1}{\int\frac{1}{\lambda(pa(x)+(1-p)b(x))+(1-\lambda
)(qa(x)+(1-q)b(x))}dF(x)}=h(\lambda p+(1-\lambda)q).
\end{align*}

\ \hfill$\square$

\bigskip

\noindent\textbf{Lemma D.11.} $f(p)$\textit{ is continuous on} $[0,$ $1]$.

\noindent\textit{Proof}. As $f(p)$ is concave, it is continuous with respect
to $0<p<1$ (See [\textbf{7}, Theorem 10.3]). It is sufficient to prove the
assertion in the case where $p\rightarrow1^{-}$. As the function $(1-p)/p$ is
strictly decreasing with respect to $p\in(0,1)$, using Lebesgue's theorem, we
obtain
\begin{align*}
\lim_{p\rightarrow1^{-}}f(p)  &  =\frac{1}{e^{r}}\lim_{p\rightarrow1^{-}}%
p\exp(\int\log(a(x)+\frac{1-p}{p}b(x))dF(x))\\
&  =\frac{1}{e^{r}}\exp(\int\log a(x)dF(x))=f(1),
\end{align*}
where $b(x)\geq0$.\hfill$\square$

\bigskip

\noindent\textbf{Lemma D.12.} $h(p)$\textit{ is continuous on} $[0,$ $1]$.

\noindent\textit{Proof}. It is not difficult to verify that $0\leq
h(p)<pE^{A}+(1-p)E^{B}<\infty$. Similar to case of Lemma D.11, we obtain the
conclusion. \hfill$\square$

\bigskip

\noindent\textbf{Lemma D.13.} \textit{Let }$\alpha(x,y)$\textit{ be continuous
with respect to }$(x,y)\in(-\delta,\delta)^{n}\times(0,\delta)$\textit{ for
some positive number }$\delta>0$\textit{, and nondecreasing with respect to
}$y\in(0,\delta)$\textit{ for each }$x\in(-\delta,\delta)^{n}$\textit{. Let
}$\beta(x)$\textit{ be continuous with respect to }$x\in(-\delta,\delta)^{n}%
$\textit{, satisfying }$\lim_{y\rightarrow0^{+}}\alpha(x,y)$\textit{ }%
$=\beta(x)$\textit{ for each }$x\in(-\delta,\delta)^{n}$\textit{. Then,
}$\alpha(x,y)$\textit{ has a unique continuous extension on }$(-\delta
,\delta)^{n}\times\lbrack0,\delta)$\textit{.}

\noindent\textit{Proof}. Define $\alpha(x,0):=\beta(x)$ $(x\in(-\delta
,\delta)^{n})$. It is sufficient to prove that $\alpha(x,y)$ is continuous at
$((0),0)$. Choose $\varepsilon>0$. (1) As $\beta(x)$ is continuous at $x=(0)$,
$0<\delta_{1}<\delta$ exists such that $\left\vert \beta(x)-\beta
((0))\right\vert <\varepsilon$ if $x\in(-\delta_{1},\delta_{1})^{n}.$ (2) As
$\lim_{y\rightarrow0^{+}}\alpha((0),y)=\beta((0)),$ $0<\delta_{2}<\delta_{1}$
exists such that $\left\vert \alpha((0),y)-\beta((0))\right\vert <\varepsilon$
if $0<y<\delta_{2}.$ (3) As $\alpha(x,y)$ is continuous at $((0),$ $\delta
_{2}/2)$, $0<\delta_{3}<\delta_{2}/2$ exists such that $\left\vert
\alpha(x,y)-\alpha((0),\delta_{2}/2)\right\vert <\varepsilon$ if $x\in
(-\delta_{3},\delta_{3})^{n}$ and $\left\vert y-\delta_{2}/2\right\vert
<\delta_{3}.$ Therefore, for each $(x^{\prime},y^{\prime})$ such that
$x^{\prime}\in(-\delta_{3},$ $\delta_{3})^{n}$ and $0\leq y^{\prime}%
<\delta_{2}/2$, we have $\beta(x^{\prime})\leq\alpha(x^{\prime},y^{\prime
})\leq\alpha(x^{\prime},\delta_{2}/2)$. Thus, $\alpha(x^{\prime},y^{\prime})$
$-\alpha((0),0)$ $\geq\beta(x^{\prime})-\beta((0))>-\varepsilon$. Moreover,
$\alpha(x^{\prime},y^{\prime})-\alpha((0),0)\leq\alpha(x^{\prime},\delta
_{2}/2)$ $-\beta((0))$ $=\alpha(x^{\prime},\delta_{2}/2)$ $-\alpha
((0),\delta_{2}/2)+\alpha((0),\delta_{2}/2)-\beta((0))<2\varepsilon$. Hence,
we have the conclusion.\hfill$\square$

\bigskip

\noindent\textbf{Remark to Lemma D.13.} Lemma D.13 is valid if the condition
$(-\delta,\delta)^{n}\times(0,\delta)$ is replaced by $[0,\delta)^{n}%
\times(0,\delta)$ and/or the term \textquotedblleft
nondecreasing\textquotedblright\ is replaced by \textquotedblleft
nondecreasing or nonincreasing.\textquotedblright

\bigskip

\noindent\textbf{Lemma D.14. }\textit{Assume }$f(p)\geq h(p)$ \textit{for each
}$p\in\lbrack0,$\textit{ }$1]$ \textit{and choose }$L:=\sup_{0<p<1}%
g(p)+1$\textit{. Then the function }%
\[
V_{t}(p,t,u):=\int\frac{pa(x)+(1-p)b(x)-u}{\left(  pa(x)+(1-p)b(x)\right)
t-ut+u}dF(x)
\]
\textit{ is upper and lower semicontinuous on }$\overline{D}=\{(p,t,u)$%
\textit{ }:\textit{ }$0\leq p\leq1$\textit{, }$0\leq t$ $\leq1$\textit{, and
}$f(p)\leq u$ $\leq L\}$\textit{. Moreover, }$V_{t}(p,t,u)=-\infty$\textit{ if
and only if }$t=1$ and $h(p)=0$\textit{.}

\noindent\textit{Proof}. As\textit{ }$f(p)\geq h(p),$ $g(p)$ exists such that
$g(p)\leq pE^{A}+(1-p)E^{B}$. Put $c(x):=pa(x)$ $+(1-p)b(x)$, then
$V_{t}(p,t,u)=\int(c(x)-u)/(c(x)t-ut+u)dF(x)$. From Hartogs' theorem,
$V_{t}(p,t,u)$ is analytic in $D$ $:=\{(p,t,u)$ : $0<p<1$, $0<t<1$, and $f(p)$
$<u<L\}$ (see [\textbf{3}, Lemma 3.1]).

First, assume that $h(p)>0$ for each $p\in\lbrack0,1].$ Put $U(x,$ $p,$ $t,$
$u):=(c(x)-u)$ $/(c(x)t$ $-ut+u)$, then $V_{t}(p,t,u)=\int U(x,p,t,u)$ $dF(x)$.

(1) From%
\[
\frac{\partial U}{\partial t}(x,p,t,u)=-\left(  \frac{c(x)-u}{c(x)t-ut+u}%
\right)  ^{2}\leq0,
\]
we can use Lebesgue's monotone theorem to obtain%
\begin{align*}
\lim_{t\rightarrow0^{+}}V_{t}(p,t,u)  &  =V_{t}(p,0,u)=\frac{pE^{A}%
+(1-p)E^{B}}{u}-1,\\
\lim_{t\rightarrow1^{-}}V_{t}(p,t,u)  &  =V_{t}(p,1,u)=1-\frac{u}{h(p)}.
\end{align*}
Notice that $V_{t}(p,0,u)$ and $V_{t}(p,1,u)$ are analytic in $\{(p,u)$ :
$0<p<1$ and $f(p)<u<L\}$ (see [\textbf{3}, Lemma 3.1]). Moreover, as $f(p)>0$
and $h(p)>0$, we have the following properties:%
\[%
\begin{array}
[c]{cc}%
\lim_{p\rightarrow0^{+}}V_{t}(p,0,u)=V_{t}(0,0,u), & \lim_{p\rightarrow0^{+}%
}V_{t}(p,1,u)=V_{t}(0,1,u),\\
\lim_{p\rightarrow1^{-}}V_{t}(p,0,u)=V_{t}(1,0,u), & \lim_{p\rightarrow1^{-}%
}V_{t}(p,1,u)=V_{t}(1,1,u),\\
\lim_{u\rightarrow f(p)^{+}}V_{t}(p,0,u)=V_{t}(p,0,f(p)), & \lim_{u\rightarrow
f(p)^{+}}V_{t}(p,1,u)=V_{t}(p,1,f(p)),\\
\lim_{u\rightarrow L^{-}}V_{t}(p,0,u)=V_{t}(p,0,L)), & \lim_{u\rightarrow
L^{-}}V_{t}(p,1,u)=V_{t}(p,1,L)).
\end{array}
\]
To obtain these equalities, we have used the inequalities $\partial
^{2}U(x,p,t,u)/\partial p^{2}\leq0$ and $\partial U(x,p,t,u)$ $/\partial
u\leq0$, which will be shown in (2) and (3).

\bigskip

(2) From%
\[
\frac{\partial^{2}U}{\partial p^{2}}(x,p,t,u)=-\frac{2ut\left(
a(x)-b(x)\right)  ^{2}}{\left(  c(x)t-ut+u\right)  ^{3}}\leq0,
\]
$U(x,p,t,u)$ is concave with respect to $p\in(0,1)$. Therefore, $U(x,p,t,u)$
is nondecreasing or nonincreasing on $(0$, $\varepsilon)$ for some positive
$\varepsilon$. Therefore, using Lebesgue's monotone theorem, we obtain
\begin{align*}
\lim_{p\rightarrow0^{+}}V_{t}(p,t,u)  &  =V_{t}(0,t,u)=\int\frac
{b(x)-u}{b(x)t-ut+u}dF(x),\\
\lim_{p\rightarrow1^{-}}V_{t}(p,t,u)  &  =V_{t}(1,t,u)=\int\frac
{a(x)-u}{a(x)t-ut+u}dF(x).
\end{align*}
Notice that $V_{t}(0,t,u)$ and $V_{t}(1,t,u)$ are analytic in $\{(t,u)$ :
$0<t<1$ and $f(p)<u<L\}$ (see [\textbf{3}, Lemma 3.1]). Moreover, as $f(p)>0$
and $h(p)>0$, we have the following properties:%
\[%
\begin{array}
[c]{cc}%
\lim_{t\rightarrow0^{+}}V_{t}(0,t,u)=V_{t}(0,0,u), & \lim_{t\rightarrow0^{+}%
}V_{t}(1,t,u)=V_{t}(1,0,u),\\
\lim_{t\rightarrow1^{-}}V_{t}(0,t,u)=V_{t}(0,1,u), & \lim_{t\rightarrow1^{-}%
}V_{t}(1,t,u)=V_{t}(1,1,u),\\
\lim_{u\rightarrow f(0)^{-}}V_{t}(0,t,u)=V_{t}(0,t,f(0)), & \lim_{u\rightarrow
f(p)^{-}}V_{t}(1,t,u)=V_{t}(1,t,f(1)),\\
\lim_{u\rightarrow L^{+}}V_{t}(0,t,u)=V_{t}(0,t,L), & \lim_{u\rightarrow
L^{+}}V_{t}(1,t,u)=V_{t}(1,t,L).
\end{array}
\]
To obtain these equalities, we have used the inequality $\partial
U(x,p,t,u)/\partial u\leq0$, which will be shown in (3).

\bigskip

(3) By the inequality%
\[
\frac{\partial U}{\partial u}(x,p,t,u)=-\frac{c(x)}{\left(  c(x)t-ut+u\right)
^{2}}\leq0,
\]
we can use Lebesgue's monotone theorem to obtain%
\begin{align*}
\lim_{u\rightarrow f(p)^{+}}V_{t}(p,t,u)  &  =V_{t}(p,t,f(p)),\\
\lim_{u\rightarrow L^{-}}V_{t}(p,t,u)  &  =V_{t}(p,t,L).
\end{align*}
As $f(p)$ is analytic in $(0,1),$ $V_{t}(p,t,f(p))$ and $V_{t}(p,t,L)$ are
analytic in $\{(p,t)$ : $0<p<1$ and $0<t<1\}$ (see [\textbf{3}, Lemma 3.1]).
Moreover, as $f(p)>0$ and $h(p)>0$, we have the following properties:%
\begin{align*}
\lim_{t\rightarrow0^{-}}V_{t}(p,t,f(p))  &  =V_{t}(p,0,f(p))=\frac
{pE^{A}+(1-p)E^{B}}{f(p)}-1,\\
\lim_{t\rightarrow1^{+}}V_{t}(p,t,f(p))  &  =V_{t}(p,1,f(p))=1-\frac
{f(p)}{h(p)},\\
\lim_{p\rightarrow0^{+}}V_{t}(p,t,L)  &  =V_{t}(0,t,L),\\
\lim_{p\rightarrow1^{-}}V_{t}(p,t,L)  &  =V_{t}(1,t,L),\\
\lim_{t\rightarrow0^{-}}V_{t}(p,t,L)  &  =V_{t}(p,0,L)=\frac{pE^{A}%
+(1-p)E^{B}}{L}-1,\\
\lim_{t\rightarrow1^{+}}V_{t}(p,t,L)  &  =V_{t}(p,1,L)=1-\frac{L}{h(p)}.
\end{align*}
By the relations
\begin{align*}
V_{t}(p,t,f(p))  &  =\frac{1}{t}-\frac{1}{t}\int\frac{1}{\frac{c(x)}%
{f(p)}t+1-t}dF(x),\\
\left\vert \frac{1}{\frac{c(x)}{f(p)}t+1-t}\right\vert  &  \leq\left\vert
\frac{1}{1-t}\right\vert ,
\end{align*}
and Lebesgue's dominated theorem, we have%
\begin{align*}
\lim_{p\rightarrow0^{+}}V_{t}(p,t,f(p))  &  =V_{t}(0,t,f(0)),\\
\lim_{p\rightarrow1^{-}}V_{t}(p,t,f(p))  &  =V_{t}(1,t,f(1)).
\end{align*}
However, these convergences are not necessarily monotonic. The continuity of
$V_{t}(p,$ $t,$ $f(p))$ near the boundaries $\{(p,t,u)$ : $p=0$, $0\leq
t\leq1$, and $u=f(0)\}$ and $\{(p,t,u)$ : $p=1$, $0\leq t\leq1$, and
$u=f(1)\}$ can be deduced, if we consider the case where $f(p)$ is replaced by
$f(p)-\min_{0\leq p\leq1}f(p)/2$.

From (1), (2), (3), Lemma D.13, and Remark to Lemma D.13, we obtain that
$V_{t}(p,t,u)$ is continuous on $\overline{D}=\{(p,t,u)$ : $0\leq p\leq1$,
$0\leq t\leq1$, and $f(p)\leq u\leq L\}.$

Second, it is easy to see that $V_{t}(p,t,u)=-\infty$ if and only if
$(p,t,u)\in M:=\{t$ $=1$ and $h(p)=0\}$. As $M$ is compact, $V_{t}(p,t,u)$ is
continuous on the open set $M^{c}$. Therefore, for each real number $m$,
$\{(p,t,u):V_{t}(p,t,u)>m\}\subset M^{c}$ is open in $\overline{D}$, which
implies that $V_{t}(p,t,u)$ is lower semicontinuous.

After this, we will prove that $V_{t}(p,t,u)$ is upper semicontinuous.

Consider the case where $h(p)=0$ $(0\leq p\leq1)$, and put $L^{\prime}%
:=\min(f(0),f(1))$. As $\partial U(x,p,t,u)/\partial u\leq0$, we obtain%
\[
V_{t}(p,t,u)=\int\frac{c(x)-u}{c(x)t-ut+u}dF(x)\leq\int\frac{c(x)-L^{\prime}%
}{c(x)t-L^{\prime}t+L^{\prime}}dF(x).
\]
Put $W(p,t):=\int\left(  c(x)-L^{\prime}\right)  /\left(  c(x)t-L^{\prime
}t+L^{\prime}\right)  dF(x)$, then, as $\partial U(x,p,t,L^{\prime})/\partial
t\leq0$, $W(p,t)$ is nonincreasing with respect to $t\in(0,1)$ and
$\lim_{t\rightarrow1^{-}}W(p,t)$ $=-\infty$. Thus, for each real number $m$
and $p\in\lbrack0,1]$, $t_{p}>0$ exists such that $1-t_{p}<t^{\prime}<1$
implies $W(p,t^{\prime})$ $<m$ $-2$. As $W(p,t)$ is continuous at $(p,$
$1-t_{p}/2)$, $0<\delta_{p}<t_{p}/2$ exists such that the conditions
$p-\delta_{p}<p^{\prime}$ $<p+\delta_{p}$, $0\leq p\leq1,$ and $1-t_{p}%
/2-\delta_{p}<t^{\prime}<1-t_{p}/2+\delta_{p}$ imply $W(p,$ $1-t_{p}/2)-1$
$<W(p^{\prime},t^{\prime})$ $<W(p,1-t_{p}/2)$ $+1$. It should be noted that
the set of open intervals $\{(p-\delta_{p},$ $p+\delta_{p})$ $\cap
\lbrack0,1]\}_{0\leq p\leq1}$ is an open covering of the compact set $[0,1]$.
Therefore, a finite subcovering $\{(p_{i}-\delta_{pi},$ $p_{i}+\delta_{p_{i}%
})\cap\lbrack0,1]\}_{i=1,2,...,m\text{ }}$exists. Put $\delta:=\min
_{i=1,2,...,m}\delta_{p_{i}}$, then for each $0\leq p^{\prime}\leq1$ and
$1-\delta<t^{\prime\prime}<1$, we have $W(p^{\prime},t^{\prime\prime})<m-1$.
It is not difficult to see that $V_{t}(p,t,u)$ is continuous on the compact
set $K:=\{(p,t,u)$ : $0\leq p\leq1$, $0\leq t$ $\leq1-\delta$, and $f(p)\leq
u$ $\leq L\}$. Therefore, $\{(p,t,u)\in\overline{D}$ $:V_{t}(p,t,u)\geq m\}$
$\subset K\subset\overline{D}$ is compact in $K$ and also in $\overline{D}$.
This implies that $\{(p,t,u):V_{t}(p,t,u)<m\}$ is open in $\overline{D}$.

Consider the case where $h(p_{0})>0$ for some $p_{0}\in\lbrack0,1]$. As $h(p)$
is concave, the compact set $\{p$ : $h(p)=0\}$ is $\{0\}$, $\{1\}$, or
$\{0,1\}$. In each case, the upper semicontinuity can be proved in a similar
fashion as above.\hfill$\square$

\bigskip

\noindent\textbf{[3, Lemma 3.1]} \textit{The function }$w_{\beta}(z):=\int
_{I}(a(x)-\beta)/(a(x)z-z\beta+\beta)dF(x)$\textit{ is analytic with respect
to two complex variables }$z:=t+si$\textit{ and }$\beta:=u+hi$\textit{ such
that}

(a)\textit{ }$\max(\varepsilon$\textit{, }$\xi)<u<L,$

(b)\textit{ }$|h|<\varepsilon^{6}/(32(L+1)R^{2}),$

(c)\textit{ }$|z|<R,$\textit{ and }$z\notin\{z:\left\vert s\right\vert
\leq\varepsilon\}\cap\{z:t\leq\varepsilon$\textit{ or }$t\geq u/(u-\xi
)-\varepsilon\},$

\textit{where }$0<\varepsilon<\min(1/2$\textit{, }$u/(2(u-\xi)))$\textit{,
}$\max(\varepsilon$\textit{, }$\xi)<L<+\infty$\textit{, }$\max(2$, $u/(u$
$-\xi))<R<+\infty$\textit{, }$i:=\sqrt{-1}$, $\xi:=\inf_{x}a(x)$,
$\operatorname{Im}(z)=s$, \textit{and} $\operatorname{Im}(\beta)=h.$

\bigskip

\noindent\textbf{Lemma D.15.}\textit{ Under the assumption of Lemma D.14, }%
\[
V(p,t,u):=\int\log(\frac{pa(x)+(1-p)b(x)}{u}t-t+1)dF(x),
\]
\textit{is continuous on }$\overline{D}=\{(p,t,u)$\textit{ }:\textit{ }$0\leq
p\leq1$\textit{, }$0\leq t\leq1$\textit{, and }$f(p)\leq u\leq L\}$\textit{. }

\noindent\textit{Proof}. As $\partial V(p,t,u)/\partial t=V_{t}(p,t,u)$ on $D$
and $V(p,t,u)<\infty$ on $\overline{D}$, using Lemmas D.13 and D.14, we have
the conclusion.\hfill$\square$

\bigskip

\noindent\textbf{Lemma D.16.} $g(p)$\textit{ is continuous on }$[0,$\textit{
}$1]$\textit{ if }$f(p)\geq h(p)$ \textit{for each }$p\in\lbrack0,$\textit{
}$1]$.

\noindent\textit{Proof}. From $f(p)\geq h(p)$, $g(p)$ exists such
that$\ g(p)\geq f(p)\geq h(p)$ (see Remark 3.1 and Lemma D.6). As $g(p)$ is
concave, it is continuous with respect to $0<p<1$ (See [\textbf{7}, Theorem
10.3]), and so $c:=\lim_{p\rightarrow1^{-}}$ $g(p)\geq g(1)$ exists. It is
sufficient to show that $c=g(1)$.

By Lemmas D.14 and D.15, the set $K:=\{(p,t,u)\in$ $\overline{D}$ :
$V(p,t,u)=r$ and $V_{t}(p,t,u)$ $=0\}$ is compact. As $c=\lim_{p\rightarrow
1^{-}}$ $g(p)$, a strictly increasing sequence $\{p_{n}\}$ exists such that
$\lim_{n\rightarrow\infty}p_{n}=1$ and $\lim_{n\rightarrow\infty}$
$g(p_{n})=c$. As the sequence $\{(p_{n},t_{g(p_{n})},$ $g(p_{n}))\}$ is in the
compact set $K$, a subsequence $\{p_{n}^{\prime}\}\subset\{p_{n}\}$ exists
such that $t_{\ast}:=\lim_{n\rightarrow\infty}t_{g(p_{n}^{\prime})}$ and $(1,$
$t_{\ast},$ $c)\in K$. By the uniqueness of the solutions of the simultaneous
equations (see Remark 3.1 and [\textbf{3}, Section 6]), we obtain $c=g(1)$.

\ \hfill$\square$

\bigskip

\noindent\textbf{Lemma D.17.} $u(p)$ \textit{is concave on} $[0,$ $1]$.

\noindent\textit{Proof}. Let $\{p,$ $q,$ $\lambda\}\subset\lbrack0,1]$, $p<q$,
and $r:=\lambda p+(1$ $-\lambda)q$. We will show that $\lambda u(p)$
$+(1-\lambda)u(q)\leq u(r)$. From Lemmas D.6 and D.9, we have $\lambda
u(p)+(1$ $-\lambda)u(q)\leq\lambda g(p)+(1-\lambda)g(q)\leq g(r)$. Therefore,
if $u(r)=g(r)$, then the assertion is proved. Henceforth, we assume that
$u(r)=f(r)<g(r)$.

In the case where $u(p)$ $=f(p)$ and $u(q)=f(q)$, the assertion follows from
Lemma D.8. In the case where $u(p)=g(p)>f(p)$ and $u(q)=g(q)>f(q)$, there are
$p<z<r$ and $r<w<q$ such that $u(z)=f(z)=g(z)=h(z)$ and $u(w)=f(w)$
$=g(w)=h(w)$. As $g(p)$ is concave, from $p<z<r<w<q$, we have%
\begin{align*}
\lambda g(p)+(1-\lambda)g(q)  &  \leq\frac{w-r}{w-z}g(z)+\frac{r-z}{w-z}g(w)\\
&  =\frac{w-r}{w-z}f(z)+\frac{r-z}{w-z}f(w)\leq f(r),
\end{align*}
which implies that $\lambda u(p)+(1-\lambda)u(q)\leq u(r)$.

The other cases can be obtained in a similar fashion.\hfill$\square$

\bigskip

\noindent\textbf{Lemma D.18.} $u(p)$ \textit{is continuous with respect to
}$p\in\lbrack0,$ $1]$.

\noindent\textit{Proof}. As $u(p)$ is concave, it is continuous with respect
to $0<p<1$ (See [\textbf{7}, Theorem 10.3]). We will only show that
$\lim_{p\rightarrow1^{-}}u(p)=u(1)$. In the case where $h(1)>f(1)$, $u(p)$
$=f(p)$ in a neighborhood of $1$. Therefore, the continuity of $f(p)$ deduces
$\lim_{p\rightarrow1^{-}}u(p)$ $=\lim_{p\rightarrow1^{-}}f(p)$ $=f(1)=u(1)$.
Similarly, in the case where $h(1)<f(1)$, we have $\lim_{p\rightarrow1^{-}%
}u(p)=\lim_{p\rightarrow1^{-}}g(p)=g(1)=u(1)$. In the case where $h(1)=f(1)$,
using Lemma D.7, we have $h(1)=\lim\inf_{p\rightarrow1^{-}}\min(f(p),$
$g(p))\leq\lim\inf_{p\rightarrow1^{-}}u(p)$ $\leq\lim\sup_{p\rightarrow1^{-}%
}u(p)$ $\leq\lim\sup_{p\rightarrow1^{-}}\max(f(p),$ $g(p))=h(1)$, which
implies the conclusion.\hfill$\square$

\bigskip

\noindent\textbf{Theorem D.19.} $u_{r}^{\sum_{i=1}^{n}p_{i}A_{i}}$ \textit{is
continuous with respect to} $(p_{i})\in Q$.

\noindent\textit{Proof}. As $u_{r}^{A}$ is finite and concave on $Q$ (Lemma
D.17), $u_{r}^{A}$ is continuous on the relative interior of $Q$ (see
[\textbf{7}, Theorem 10.1]) and has a unique continuous extension on $Q$ (see
[\textbf{7}, Theorems 10.3 and 20.5]). Therefore, we need to show the relation
$\lim_{p\rightarrow1^{-}}u_{r}^{pA+(1-p)B}=u_{r}^{A}$, where $A$ or $B$ is a
relative boundary point or a relative interior point of $Q$, respectively. In
this instance, Lemma D.18 leads to the conclusion.\hfill$\square$

\bigskip

\noindent\textbf{[3, Example 6.6]} The European put option is given by%
\[
a(x)=\max(K-Se^{rT}e^{x},0),\text{ \ \ }dF(x)=\frac{1}{\sqrt{2\pi T}\sigma
}e^{-\frac{(x+\sigma^{2}T/2)^{2}}{2\sigma^{2}T}}dx.
\]
We assume that the stock price $Y=Se^{rT}e^{X}$ is lognormally distributed
with volatility $\sigma\sqrt{T}$, where $S$ is the current stock price, $r$ is
the continuously compounded interest rate, $K$ is the exercise price of the
put option, and $T$ is the exercise period. The expectation $E$ is given by
\begin{align*}
E  &  =\frac{1}{\sqrt{2\pi T}\sigma}\int_{-\infty}^{\log\frac{K}{S}%
-rT}(K-Se^{rT}e^{x})e^{-\frac{(x+\sigma^{2}T/2)^{2}}{2\sigma^{2}T}}dx\\
&  =KN\left(  -\frac{\log\frac{S}{K}+(r-\frac{\sigma^{2}}{2})T}{\sigma\sqrt
{T}}\right)  -Se^{rT}N\left(  -\frac{\log\frac{S}{K}+(r+\frac{\sigma^{2}}%
{2})T}{\sigma\sqrt{T}}\right)  ,
\end{align*}
where $N(x)=\int_{-\infty}^{x}e^{-x^{2}/2}/\sqrt{2\pi}dx$ is the cumulative
standard normal distribution function.

When $S=90$, $K=120$, $T=2$, $\sigma=0.1$, and $r=0.04$, we have $\xi=0$, $E$
$\fallingdotseq22.9848,$ and $H:=\int1/a(x)dF(x)=+\infty$. Therefore, from
Theorems 4.1 and 5.1, $G_{u}(t_{u})$\ ($u\in(0,$ $E)$) strictly decreases from
$+\infty$ to $1$. The equations $w_{u}(t_{u})=0$ and $G_{u}(t_{u})=e^{0.08}$
yield the price $u\fallingdotseq17.8157$. With this price, if investors
continue to invest $t_{u}$ $\fallingdotseq0.5434$\ of their current capital,
they can maximize the limit expectation of growth rate to $e^{0.08}%
\fallingdotseq1.0833.$

In general, the equation $E/u=e^{rT}$ yields the price%
\[
u=Ke^{-rT}N\left(  -\frac{\log\frac{S}{K}+(r-\frac{\sigma^{2}}{2})T}%
{\sigma\sqrt{T}}\right)  -SN\left(  -\frac{\log\frac{S}{K}+(r+\frac{\sigma
^{2}}{2})T}{\sigma\sqrt{T}}\right)  ,
\]
which is the Black-Scholes formula for the European put option. Substituting
the above mentioned values for this formula, we obtain the (higher) price
$u\fallingdotseq21.2176$ $(>17.8157)$. With this price, if the investors
continue to invest $t_{u}\fallingdotseq0.2278$ of their current capital, they
can maximize the limit expectation of growth rate to $1.0096$ $(<1.0833).$

\bigskip

\bigskip

\bigskip

\begin{center}
\textbf{References}
\end{center}

\noindent\lbrack1] \ F. Black and M. Scholes, The pricing of options and
corporate liabilities, J.

Political Economy 81 (1973), 637--654.

\noindent\lbrack2] \ G. H. Hardy, J. E. Littlewood and G. P\'{o}lya,
Inequalities, Cambridge University

Press, Reprinted 1973.

\noindent\lbrack3] \ Y. Hirashita, Game pricing and double sequence of random variables,

Preprint, arXiv:math.OC/0703076 (2007).

\noindent\lbrack4] \ Y. Hirashita, Least-Squares Prices of Games,
Preprint,\ arXiv:math.OC/0703079 

(2007).

\noindent\lbrack5] \ J. L. Kelly, Jr., A new interpretation of information
rate,\ Bell. System Tech.

J. 35 (1956), 917--926.

\noindent\lbrack6] \ D. G. Luenberger, Investment science, Oxford University
Press, Oxford, 1998.

\noindent\lbrack7] \ R. T. Rockafellar, Convex analysis, Princeton University
Press, Princeton, 1970.

\noindent\lbrack8]\ \ I. E. Schochetman, Pointwise versions of the maximum
theorem with applications

in optimization, Appl. Math. Lett. 3 (1990), 89--92.

\bigskip

\bigskip

\noindent Chukyo University, Nagoya 466-8666, Japan

\noindent yukioh@cnc.chukyo-u.ac.jp
\end{document}